\documentclass[12pt,bezier]{article}
\usepackage{amsmath,amsthm,amssymb,amsfonts,euscript,graphicx}

\theoremstyle{plain}
\newtheorem{theorem}{Theorem}

\newtheorem{definition}{Definition}

\newtheorem{conjecture}{Conjecture}
\newtheorem{problem}{Problem}

\title{\large \bf A Note on Matching in Groups and Field Extensions
}

\author{{\normalsize{\sc Babak Hassanzadeh}}\\
{\footnotesize{\it Department of Mathematics,  Azarbaijan shahid Madani University,  }}\\[-3mm]
{\footnotesize{\it Tbriz, Iran}}\\
{\footnotesize{E-mail address: $\mathsf{babakmath777@gmail.com}$}}\\
}

\date{}
\begin{document}

\maketitle
\begin{abstract}
\noindent
 The purpose of this note is to give a number of  open problems on matching theory and their relation to the well-known results in this area. We also give a linear analogue of the acyclic matchings.
\end{abstract}
\vspace{9mm}
\section{Results and Discussion}
\begin{definition}
Assume that  $G$ is an abelian group and $A$ and $B$ are two non-empty finite subsets of $G$ with $|A|=|B|$. A matching from $A$ to $B$ is a bijection $f:A\to B$ satisfying $a+f(a)\not\in A$, for any $a\in A$. Also for any matching $f:A\to B$, we define the function $m_f$ as follows:
\[\begin{array}{rl} m_f: G&\to \mathbb{N}\cup\{0\}\\ x&\mapsto \left|\left\{a\in A;\; a+f(a)=x\right\}\right|.
\end{array}\]
A matching $f$ is called to be acyclic if for any matching $g:A\to B$ with $m_f=m_g$, we have $f=g$.\\
We say $G$ has the (acyclic) matching property if for any finite non-empty subsets $A$ and $B$ of $G$ with $|A|=|B|$ and $0\not\in B$, there exists at least one (acyclic) matching from $A$ to $B$.
\end{definition}
We have the following theorems about the  matching property:
\begin{theorem}(\cite[Theorem 3.1]{8})
An abelian  group $G$ has the matching property if and only if $G$ is either torsion free or cyclic of prime order.
\end{theorem}
\begin{theorem}(\cite[Theorem 4.1]{8})
If $G$ is abelian torsion-free group, then $G$  has the acyclic matching property.
\end{theorem}
\begin{theorem}(\cite[Prop 2.2 and 2.3]{2})
There are infinitely many prime $p$ for which $\mathbb{Z}/p\mathbb{Z}$ has no acyclic matching property.
\end{theorem}
\begin{conjecture}
There are infinitely many prime $p$ for which $\mathbb{Z}/p\mathbb{Z}$ has the acyclic matching property.
\end{conjecture}
\begin{problem}
It would be interesting to find some families of primes $\{p_i\}$ such that $\mathbb{Z}/p_i\mathbb{Z}$ has no matching property. (We already found two such families in \cite{2}.) One can use the local matching property in groups to formulate such primes \cite{3}.
\end{problem}
\begin{definition}
Let $K\subseteq L$ be a field extension and $A$ and $B$ be two $n$-dimensional $K$-subspaces of $L$ and $\mathcal{A}=\{a_2,\ldots,a_n\}$ and $\mathcal{B}=\{b_1,\ldots,b_n\}$ be bases for $A$ and $B$, respectively. We say $\mathcal{A}$ is matched to $\mathcal{B}$ if 
\[a_i^{-1}A\cap B\subseteq \langle b_1,\ldots,\hat{b_i},\ldots,b_n\rangle,\]
for any $1\leq i\leq n$.\\
Also, we say $A$ is matched to $B$ if for any basis $\mathcal{A}$ of $A$, there exists a basis $\mathcal{B}$ of $B$ such that $\mathcal{A}$ is matched to $\mathcal{B}$.
\end{definition}
\begin{definition}
A strong matching from $A$ to $B$ is a linear isomorphism $\varphi:A\to B$ such that any basis $\mathcal{A}$ of $A$ is matched to the  basis $\varphi(\mathcal{A})$ of $B$. \end{definition}
We have the following theorems:
\begin{theorem}(\cite[Theorem 6.3]{7})
Let $A$ and $B$ be $n$-dimensional $K$-subspaces of $L$ distinct from $\{0\}$. There is a strong matching from $A$ to $B$ if and only if $AB\cap A=\{0\}$. In this case, any isomorphism $\varphi:A\to B$ is a strong matching.
\end{theorem}
\begin{definition}
If $K\subseteq L$ is a field extension, we say that it has the linear matching  property if for any two $n$-dimensional $K$-subspaces $A$ and $B$ of $L$ with $n\geq1$ and $1\not\in B$, $A$ is  matched to $B$.
\end{definition}
The following  theorem is about the characterization of all field extensions satisfying the linear matching property.
\begin{theorem}(\cite[Theorem 5.2]{7})
Let $K\subseteq L$ be a field extension. Then it has the linear matching property if and only if $L$ contains no proper finite-dimensional extension over $K$. See also \cite{1} for the more precise version of THeorem 5.
\end{theorem}
Our definition of the linear acyclic matching property:
\begin{definition}
We say that two linear isomorphism $f,g:A\to B$ are equivalent if there exists a linear automorphism $\phi:A\to A$ such that for all $a\in A$ one has $af(a)=\phi(a)g(\phi(a))$, and two strong matchings $f,g:A\to B$ are equivalent if they are equivalent as linear isomorphism. Then, we define a acyclic matching from $A$ to $B$ to be a strong matching $f:A\to B$ such that for any strong matching $g:A\to B$ that is equivalent to $f$, one has $f=cg$, for some constant $c\in K$. Finally we say that $K\subseteq L$ satisfies the linear acyclic matching property if for every pair $A$ nad $B$ of non-zero equi-dimensional $K$-subspace of $L$ with $AB\cap A=\{0\}$, there is at least one acyclic matching from $A$ to $B$. See [4 and 5] for more details on linear analogous of matchings.
\end{definition}
\begin{theorem}(\cite[Theorem 4.5]{1})
Let $K\subseteq L$ be a purely transcendental field extension of a field $K$. Then $K\subseteq L$ satisfies the linear acyclic matching property.
\end{theorem}
\begin{conjecture}
There are infinitely many prime $p$ such that there exists a field extension $K\subseteq L$ with $[L:K]=p$ and no linear acyclic matching property.
\end{conjecture}
The main possible key to solve Conjecture 2 is an elementary linear algebra theorem which states that a finite dimensional vector space over an infinite field cannot be written as a finite union of its proper subspaces \cite{6}.

\end{document}